# MINIMAL PERIOD OF SOLUTIONS TO LIPSCHITZ DIFFERENTIAL EQUATINS WITH ARBITRARY VECTOR NORM

A.A. Zevin

*Institute of Transport System and Technologies,*

*National Academy of Sciences of Ukraine*

*Pisarzhevsky 5. 49005, Dnepr, Ukraine*

e-mail: zevin@westa-inter.com

**Abstract**

The Lipschitz differential equation, $\dot{x} = f(x)$, in spaces $X \in \mathrm{C}^n$ and $X \in \mathrm{R}^n$ is considered. The minimal period problem is to find the exact lower bound for periods of non-constant solutions, expressed in the Lipschitz constant $L$. In this paper, some inequality for the components, $x_k(t)$, which is independent on the space $X$, is found. As a result, it is proved that for any $X$ and $n$, the normalized minimal period, $k = TL \geq 2\pi$. In the space $\mathrm{C}^n$, the equality $k = 2\pi$ is reached for any $X$. For $X \in \mathrm{R}^n$, this equality is attained for univercally adopted norms.

## 1. Introduction and preliminary result

We consider periodic solutions $x(t) = x(t+T) \neq const$ of the equation
$$\dot{x} = f(x) \qquad (1.1)$$
in spaces $X \in \mathrm{C}^n$ and $X \in \mathrm{R}^n$. The function $f(x)$ satisfies the Lipschitz condition,
$$\|f(x') - f(x'')\| \leq L \|x' - x''\| \qquad (1.2)$$

The normalized minimal period, $k = TL$. of solutions $x(t)$ depends only on the space $X$.

The first exact solution of the minimal period problem belongs to Yorke [1], who proved that in the space $R^n$ with the Euclidean norm, $k = 2\pi$. The equality is achieved at the system

$$\dot{x}_1 = Lx_2, \quad \dot{x}_2 = -Lx_1 \tag{1.3}$$

Subsequently, Lasota and Yorke [2] generalized this result to any Gilbert space.

The equality $k = 2\pi$ is also valid (Mawhin and Walter, [3]) for equations of a higher order,

$$x^{(r)} = f(x), x \in R^n \tag{1.4}$$

(note that here $T = k/L^{1/r}$).

Busenberg et all Fisher and Martelli found [4] that in the general Banach space, $k = 6$. For the family of $l^p$-norms, Nieuwenhuis, Robinson and Steineberger found bounds for $k$ [6] that are strictly larger than 6. For equation (1.4) with even $r$, the above question was completely answered by Zevin, [7]) where it was proved that $k = 2\pi$ for any vector norm. The same result was claimed for equation (1.1) [7,8], however, the corresponding proofs are incorrect.

Below it is proved that the equality $k = 2\pi$ is attained in any space $X \in C^n$. For $X \in R^n$, the spaces are indicated for which this bound is also achieved (some of these results were recently obtained in [9] using different arguments). In the further proofs, the following lemma plays a key role.

Let $x(t)$ be a periodic solution to system (1.1), (1.2),

$$z_k(t) = x_k(t) - x_k(t+\tau) \tag{1.5}$$

Lemma 1. The function $z_k(t)$ satisfies the inequality

$$|\dot{z}_k(t)| \leq L|z_k(t)| \tag{1.6}$$

Proof. Putting in (1.2)

$$x'_k = x_k(t), \quad x''_k = x_k(t+\tau),$$

from (1.1) and (1.5) we have

$$\|f(x') - f(x'')\| = \|\dot{z}(t)\| \leq L\|z(t)\|, \tag{1.7}$$

By definition, inequality (1.2) is valid for any $x'_i, x''_i$ and, therefore, for any $z$. So, putting in (1.7) $z_i = 0$ for $i \neq k$ and observing that $|\dot{z}_k(t)| \leq \|\dot{z}(t)\|$, we obtain inequality (1.6). □

2. **Main results**

Suppose first that $x \in C^n$.

Theorem 1. The minimal period of non-constant solutions to system (1.1), (1.2) in any space $X \in C^n$ equals

$$T_* = \frac{2\pi}{L}, \qquad (2.1)$$

Proof. Putting in (1.6) $z_k(t) = u(t) + iv(t)$, we obtain

$$\int_0^T |\dot{z}_k(t)|^2 dt = \int_0^T (\dot{u}(t)^2 + \dot{v}(t)^2) dt \leq L^2 \int_0^T |z_k(t)|^2 dt = L^2 \int_0^T (u(t)^2 + v(t)^2) dt \qquad (2.2)$$

As is clear from (1.5), the mean value of $z_k(t)$ and, therefore, $u(t)$ and $v(t)$ on $[0,T]$ is equal to zero. Then, as is known, for $T = 2\pi$ and $L = 1$, inequality (2.2) is strict. Therefore, for $L \neq 1$, inequality (2.2) implies $TL \geq 2\pi$.

The equality in (2.2) is achieved when $z_k(t) = c_k \exp(iLt)$. So, for any $X \in C^n$, the equality in (2.1) holds at the equation

$$\dot{x} = iLx \qquad (2.3)$$

consisting of $n$ uncoupled equations of the first order. □

Clearly, the equality in (2.1) is also valid for the equation

$$\dot{x} = Ax, \qquad (2.4)$$

if the matrix $A$ has an eigenvalue $\lambda_i = \pm iL$ with $L$ being the Lipschitz constant of the function $Ax$. As is known, $L = \|A\|$ where $\|A\|$ is an induced norm of the matrix $A$,

$$\|A\| = \sup_{z \neq 0} \frac{\|Az\|}{\|z\|}$$

Since for any metric, $\|A\| \geq |\lambda_i|$, $i = 1,..,n$, the normalized period takes the minimal value, $k = TL = 2\pi$, when

$$\|A\| = |\lambda_j| = \max |\lambda_i| \qquad (2.5)$$

If a matrix $A$ satisfies equality (2.5) but $\lambda_j = L\exp(i\mu)$ with $\mu \neq \pm\pi/2$, one can take $A(\mu) = \exp(-i\mu \pm \pi/2)A$. Clearly, $\|A(\mu)\| = \|A\|$ and the corresponding eigenvalue, $\lambda_j(\mu) = iL$, so, the normalized minimal period is attained at the equation $\dot{x} = A(\mu)x$.

Let $A = \mathrm{diag}(\lambda_1,...,\lambda_n)$, then condition (2.5) is satisfied, so, $k = 2\pi$. Clearly, the last is also true when a matrix $A$ is reduced to the diagonal form by a transformation that preserves its norm.

In particular, the normal matrix ($A^*A = AA^*$) is reduced to the diagonal form by a unitary transformation that preserves the Euclidean norm [10]. Therefore, in the space with the Euclidean norm, the minimal period, $k = 2\pi$, is attained at equation (2.4) with a normal matrix $A$.

For example, the anti-Hermitian matrix ($A^* = -A$) is normal and all its non-zero eigenvalues are imaginary [10]. Therefore, in the case of Euclidean norm, the minimal period is reached at equation (2.4) with any anti-Hermitian matrix $A$.

Let now considera real space $X \in \mathbf{R}^n$ where $n > 1$, because for $n = 1$, non-constant periodic solutions do not exist. From Theorem 1 it follows that here the periods satisfy the inequality

$$T \geq \frac{2\pi}{L}. \tag{2.6}$$

The problem is to find spaces for which the equality in (2.6) is attained.

For a real equation (2.4), the Jordan block, corresponding to an eigenvalue $\lambda_k = \pm iL$, is of the form

$$J = \begin{bmatrix} 0 & L \\ -L & 0 \end{bmatrix} \tag{2.7}$$

Hence, the equality in (2.6) is achieved if $\|J\| = L$. Obviously, this is the case when the norm of the vector $x^1 = (x_1, x_2)$ is equal to that of $x^2 = (-x_2, x_1)$. Note that this condition is satisfied for the most common norms, e.g., ,

$$\|x\|_{l^p} = \left( \sum_{j=1}^{n} |x_j| \right)^{1/p} \quad \text{and} \quad \|x\|_{\infty} = \max_j |x_j|$$

In particular, for $p = 2$, the $l^p$-norm is Euclidean. So, analogously to the space $\mathbf{C}^n$, we find that in the space $\mathbf{R}^n$ with the Euclidean norm, the minimal period is attained at equation $\dot{x} = Ax$ with any anti-symmetric ($A' = -A$) matrix.

### 3. Conclusion

The main results of this paper are as follows.

1. For components $x_k(t)$ of a solution $x(t)$ to Lipschitz differential equation (1.1), inequality (1.6) is found. It allows to obtain estimates for solutions that are independent of the system dimention $n$ and the space $X$.

2. For $X \in \mathbf{C}^n$, it is proved that the normalized minimal period equals $2\pi$. In general, it is attained at an equation $\dot{x} = Ax$ when the matrix $A$ is reduced to a diagonal form by a transformation, preserving its norm.

3. For $X \in \mathbf{R}^n$, the normalized minimal period $k \geq 2\pi$; the equality holds if $\|x\|$ is invariant on permutation of the indices.